\documentclass[12pt,a4paper]{article}

\usepackage{graphicx}
\usepackage{amsmath,amsfonts}
\usepackage{mathrsfs}

\newtheorem{theorem}{Theorem}
\newtheorem{lemma}{Lemma}
\newtheorem{proposition}{Proposition}

\let\ds\displaystyle
\let\scr\mathscr
\let\goth\mathfrak
\let\ca\mathcal
\def\zT{{\vphantom{\widetilde T} T}}
\def\zs#1{_{\lower 3pt \hbox{$\scriptstyle#1$}}}
\def\Pb{\mathbf{P}}
\def\Ex{\mathbf{E}}
\def\RR{\mathbb{R}}
\def\taille{6cm}
\def\twopic#1#2#3#4{\nobreak\bigskip\leavevmode
\hbox to \taille{\vbox{%
\hbox {\includegraphics*[width=\taille]{#1}}%
\hbox to \taille{\hfill\small #2\hfill}}}%
\hbox to 0.5cm{\hss}
\hbox to \taille{\vbox{%
\hbox to \taille{\includegraphics*[width=\taille]{#3}}%
\hbox to \taille{\hfill\small #4\hfill}}}\bigskip}

\begin{document}
\title{Hypotheses Testing: Poisson Versus Self-Exciting}

\author{Sergue\"{\i} \textsc{Dachian}\\
{\small Laboratoire de Math\'ematiques, Universit\'e Blaise Pascal}\\[8pt]
Yury A. \textsc{Kutoyants}\\
{\small Laboratoire de Statistique et Processus, Universit\'e du Maine}}

\date{}

\maketitle
{Running title : Poisson versus self-exciting}
\begin{abstract}

We consider the problem of hypotheses testing with the basic simple
hypothesis: observed sequence of points corresponds to stationary Poisson
process with known intensity. The alternatives are stationary self-exciting
point processes. We consider one-sided parametric and one-sided nonparametric
composite alternatives and construct locally asymptotically uniformly most
powerful tests.  The results of numerical simulations of the tests are
presented.
\end{abstract}



\section{Introduction}

Let $\left\{t_1, t_2,\ldots \right\}$ be a sequence of events of a stationary
point process $X=\left\{X_t,\;t\geq 0\right\}$ ($X_t$ is a counting
process). The simplest stationary point process is, of course, Poisson process
with a constant intensity $S>0$, i.e., the increments of $X$ on disjoint
intervals are independent and distributed according to Poisson law
$$
\ds \Pb\left\{X_t-X_s=k\right\}=\frac{S^k\left(t-s\right)^k}{k!}\;
{\rm e}^{-S\left(t-s\right) },\quad 0 \leq s<t,\quad k=0,1,\ldots
$$
Therefore if we have a stationary sequence of events it is interesting to
check first of all if this model (Poisson process) corresponds well to the
observations.  The importance of this problem was discussed by Cox and Lewis
(1966), Section 6.3.

The alternatives close to the basic hypothesis correspond to the case when the
non-poissonian behavior is due to small perturbations of the Poisson process
and are the most interesting to test. For ``far alternatives'' any reasonable
test has power function close to 1 and the comparison of tests seems less
important.  Let us consider the problem of small signals detection by the
tests of fixed size $\varepsilon \in \left(0,1\right)$. Using the terminology
of statistical radiothechnics we say that there is at least two types of close
alternatives: the first one corresponds to small ``signal-noise ratio''
(signals of small energy) and the second, when the amplitude of the signal can
be small, but the total energy due to the sufficiently long time of
observation is comparable with the noise energy (see, e.g., Kutoyants
(1976)). For the first class of alternatives the approach of {\sl
locally optimal tests}, which provides the optimality of the power function at
the small vicinity of the basic hypothesis (the values of the power function
are close to $\varepsilon $) was developed (see, e.g. Capon (1961)) and
for the second class of {\sl contiguous alternatives} the optimality of the
test for a wider class of close alternatives (the values of the power function
are in $\left(\varepsilon,1 \right)$) is proved (Pitman's (1948) approach,
Le Cam's (1956) theory).

For stationary point processes with Poisson hypothesis and stationary
alternatives Davies (1977) proposed the {\sl locally optimal (efficient)}
or {\sl asymptotically locally efficient} test. This test is based on the
comparison of the derivative of the log-likelihood ratio with some threshold.
See as well Daley and Vere-Jones (2003), Section 13.1, where the
approach of Davies is discussed.

In the present note we suppose that we have observations of the point process
$X^T=\left\{X_t,\;0\leq t\leq T\right\}$ on the interval $\left[0,T\right]$
and consider two problems of hypotheses testing in the asymptotics of large
samples $\left(T\rightarrow \infty \right)$. In both problems the basic
hypothesis is simple: the observed process is standard Poisson with known
constant intensity $S_*>0$. The composite alternatives are: the observed
process is a realization of self-exciting point process (sometimes called
Hawkes (1972) process) with in the first case intensity function depending on
one-dimensional parameter and in the second case the intensity function
belonging to a wider (nonparametric) class of functions. We follow the
mentioned above Pitman-Le Cam's approach. We start with the {\sl locally
asymptotically uniformly most powerful test} (LAUMPT) in the parametric case
and the main result of the presented work is the LAUMPT where the optimality
is shown for sufficiently large class of local nonparametric alternatives.
The similar results for diffusion processes can be found in Iacus, Kutoyants
(2001) (small noise asymptotics) and Kutoyants (2003) (ergodic processes).

\section{Preliminaries}

Remind several facts from the theory of point processes (the details can be
found, for example, in Liptser and Shiryaev (2001), Chapter 18). Let $\left(\Omega ,{\goth
F},\Pb\right)$ be a probability space and let $\left\{{\goth F}_t, \;t\geq
0\right\}$ be a nondecreasing family of right continuous $\sigma $-algebras
${\goth F}_s\subset {\goth F}_t\subset {\goth F} $ for any $0\leq s<t$.  We
denote by $t_1, t_2,\ldots $ a sequence of Markov stopping times adapted to
$\left\{{\goth F}_t, \;t\geq 0\right\}$ (that means $ \left\{\omega : t_i\leq
t\right\}\in {\goth F}_t $ for all $t\geq 0$). Let $X_t$ be the number of {\sl
events $t_i$} up to time $t$, i.e., $X=\left\{X_t, {\goth F}_t,t\geq
0\right\}$ is a random process such that
$$
X_t=\sum_{i\geq 1}\chi\zs{\left\{t_i< t\right\}},\qquad t\geq 0,
$$
where $\chi\zs{\left\{A\right\}} $ is the indicator-function of the event $A$.

We assume that $\Ex X_t< \infty $ (there is no accumulation points on any
 bounded interval). The process $X$ admits a unique (up to stochastic
 equivalence) decomposition (Doob-Meyer decomposition)
 \begin{equation}
\label{1}
X_t={\cal A}_t+{\cal M}_t, \end{equation} where ${\cal M}=\left\{{\cal
 M}_t,{\goth F}_t,t\geq 0\right\}$ is a martingale and ${\cal A}=\left\{{\cal
 A}_t,{\goth F}_t,t\geq 0\right\}$ is predictable increasing process (Liptser
 and Shiryaeyv (2001), Theorem 18.1).  We suppose that the compensator ${\cal
 A} $ is absolutely continuous
$$
{\cal A}_t=\int_{0}^{t}S\left(
v,\omega \right)\;{\rm d}v,\qquad t\geq 0
$$
where ${\cal S}=\left\{S\left(t,\omega \right),{\goth F}_t,t\geq 0 \right\}$
 is called intensity function. We suppose as well that \eqref{1} is the {\sl
 minimal representation} of the point process, i.e., $S\left( t,\omega \right)$
 is measurable w.r.t. $\sigma $-algebra generated by $\{X_s,\;s< t\}$ for any
 $t>0$ and we write $S\left( t,\omega \right)=S\left( t,X\right)$. To describe
 a point process it is sufficient to specify its intensity function. We study
 in this work a special class of point processes with intensity
 functions which can be written as stochastic integrals with respect to the
 past of the underlying point process.

 In the particular case when ${\cal S}$ is deterministic, the process $X$ is
 (inhomogeneous) Poisson process with  intensity function $ S\left(v,\omega
 \right)=S\left(v\right)$.  In this case
$$
\Pb\left\{X_t-X_s=k\right\}=\frac{\left[\int_{s}^{t}S\left(v\right)
\;{\rm  d}v \right]^k}{k!}\;
 \exp{\left\{-\int_{s}^{t}S\left(v\right)\;{\rm
 d}v\right\}}
$$
for any $ t>s\geq 0$ and $ k=0,1,\ldots$. If the assumption of the
independence of increments is no more valid, then ${\cal S}$ is no more
deterministic and $X$ can be a stationary point process (see Brillinger (1975) and
Daley and Vere-Jones (2003) and references therein for wide classes of such processes and
their applications in real problems).

Remind that the distribution $\Pb_S^{\left(T\right)}$ of the point process in
the space of its realizations $\left({\scr D}\left(0,T\right), {\goth
B}_T\right)$ is entirely characterized by its intensity function ${\cal S}$.
The likelihood ratio formula (w.r.t. Poisson process of constant intensity
$S_* $) has the following form (see   Liptser and Shiriyev (2001), Theorem 19.10)
\begin{eqnarray*}
L\left(X^T\right)=\exp\left\{\int_{0}^{T}\ln \frac{S\left(t,\omega
\right)}{S_*}\;{\rm
d}X_t-\int_{0}^{T}\left[S\left(t,\omega \right)-S_*\right]{\rm d}t\right\},
\end{eqnarray*}
where we suppose that the intensity  $S\left(t,\omega\right)$ is left
 continuous  function and
$$
\Pb\left\{\int_{0}^{T}S\left(t,\omega \right){\rm
d}t<\infty \right\}
$$
under all alternatives studied in this work.

\section{One-sided parametric alternative}

Suppose that we observe a trajectory $X^T=\left\{X_t,0\leq t\leq T\right\}$ of
point process of intensity function $S_T\left(\vartheta
\right)=\left\{S\left(\vartheta ,t,\omega \right),0\leq t\leq T \right\}$. If
$\vartheta =0$, then $S\left(0 ,t,\omega \right)\equiv S_*$, i.e., this point
process is Poisson process of intensity $S_*>0$. Under alternative $\vartheta
>0$ and $S_T\left(\vartheta \right)$ is the intensity function of
self-exciting point process. As usual in such problems, we consider contiguous
alternatives (Pitman's (1948) alternatives, Roussas (1972)), hence we change
the variable $\vartheta =u/\sqrt{T}$ and test the following two hypotheses
\begin{eqnarray*}
&&{\scr   H}_0:\qquad \qquad u=0\\
 &&{\scr H}_1:\qquad \qquad u>0.
\end{eqnarray*}

We denote $\Ex_0$ the mathematical expectation under the  hypothesis ${\scr H}_0$,  and
$\Ex_u$ under (simple) alternative $\vartheta   =u/\sqrt{T}$.

Let us fix $\varepsilon \in \left(0,1\right)$ and denote by ${\scr
K}_\varepsilon $ the class of test functions $\phi_T\left(X^T\right)$ of
asymptotic size $\varepsilon $, i.e., for $\phi_T\in {\scr K}_\varepsilon $
we have
\begin{equation}
\label{2}
\lim_{T\rightarrow \infty }\Ex_0\;\phi_T\left(X^T\right)=\varepsilon.
\end{equation}
As usual, $ \phi_T\left(X^T\right)$ is the probability to accept the
hypothesis ${\scr H}_1 $ having observations
$X^T$. The corresponding power function is
$$
\beta _\zT\left(u,\phi_\zT\right)=\Ex_u\;\phi_T\left(X^T\right),\quad u\geq 0.
$$

We introduce the asymptotic optimality of  tests  with the help of the
following definition Le Cam (1956).

\bigskip
\noindent
{\bf Definition 1.} A test $\phi^*_\zT\left(\cdot \right)$ is called {\sl
locally asymptotically uniformly most powerful in the class ${\scr
K}_\varepsilon $} if for any other test $\phi_\zT\left(\cdot \right)\in {\scr
K}_\varepsilon $ and any constant $K>0$ we have
$$
\lim_{T\rightarrow \infty }\inf_{0< u\leq K}\left[\beta
_\zT\left(u,\phi_\zT^*\right)-\beta _\zT\left(u,\phi_\zT\right) \right]\geq 0 .
$$

Our goal is to  construct locally asymptotically uniformly most powerful test
 in class  ${\scr K}_\varepsilon $.

{\sl Self-exciting} type processes were introduced by Hawkes (1972) and defined by
 intensity function of the following form
\begin{equation}
\label{3}
S\left(t,\omega \right)=S_*+\int_{0}^{t-}g\left(t-s\right)\;{\rm
d}X_s=S_*+\sum_{t_i<t}^{}g\left(t-t_i\right),
\end{equation}
where  $S_*>0$, $t_i$ are the events of the point process and the function $g\left(\cdot \right)\geq 0$ satisfies
the condition
\begin{equation}
\label{4}
\rho =\int_{0}^{\infty }g\left(t\right)\;{\rm d}t<1.
\end{equation}
Remind that according to this representation of the intensity function, the
distribution of $t_1$ is exponential at rate $S_*$ and for all $n\geq 1$
$$
\Pb\left\{t_{n+1}>t\; |\;
t_1,\ldots,t_n\right\}=\exp\left(-S_*t-\int_{0}^{t}
\sum_{i=1}^{X_s}g\left(s-t_i\right)\,{\rm
d}s \right) .
$$

 Note that $\Lambda \left(t\right)=\Ex X_t$ is solution of the
 equation
\begin{eqnarray*}
\Lambda\left(t\right)&=&\Ex \int_{0}^{t} S\left(v,\omega \right)\;{\rm d}v=
S_*t+\Ex \int_{0}^{t}\int_{0 }^{v-}g\left(v-s\right)\;{\rm d}X_s\;{\rm
d}v=\\
&=&S_*t+\int_{0}^{t}\int_{0}^{v}g\left(v-s\right)\Lambda'\left(s\right)\; {\rm
d}v\;{\rm d}s.
\end{eqnarray*}
 In stationary case the intensity $S\left(t,\omega \right), $ is a stationary process
$$
S\left(t,\omega \right)=S_*+\int_{-\infty }^{t-}g\left(t-s\right)\;{\rm
d}X_s,
$$
and
$$
\Lambda \left(t\right)=\frac{ S_*}{1-\rho }\;t\equiv \mu\; t.
$$

The spectral density of this process
is
$$
f\left(\lambda \right)=\frac{\mu }{2\pi \left|1-G\left(\lambda \right)\right|^2}
$$
where
$$
G\left(\lambda \right)=\int_{0}^{\infty }{\rm e}^{{\rm i}\lambda t}g\left(t\right)\;{\rm
d}t,\qquad  \rho =G\left(0 \right).
$$

\bigskip
\noindent
{\it Example 1.} Let $g\left(t\right)=\alpha {\rm e}^{-\gamma t}$, where
$\alpha >0, \gamma >0$ and $\alpha /\gamma <1 $. Then the point process $X$
with intensity function
$$
S\left(t,\omega \right)=S_*+\alpha\; \sum_{t_i\leq t}^{}{\rm e}^{-\gamma
\left(t-t_i\right)}
$$
is self-exciting with the rate
$$
\mu =\frac{S_*\;\gamma}{\gamma -\alpha }.
$$

\bigskip
\noindent {\it Example 2.} The function $g\left(\cdot \right)$ can be chosen in such a way that the spectral
density of the point process will be rational
$$
f\left(\lambda \right)=\frac{\mu }{2\pi }\;\frac{\left|Q\left({\rm i}\lambda
\right)\right|^2 }{\left|P\left({\rm i}\lambda
\right)\right|^2}
$$
where $Q\left(z \right)=z^p+a_1z^{p-1}+\ldots+a_p$ and $P\left(z
\right)=z^p+b_1z^{p-1}+\ldots+b_p$. It is supposed that $P\left(\cdot \right)$
and $Q\left(\cdot \right)$ have no zeroes in common and no zeroes in the
closed right half plane (see Pham (1981), where the asymptotic properties
of the MLE for this model are described).

\bigskip

We assume that the observed process is either Poisson with constant intensity $S_*$ or self-exciting with {\sl
contiguous } intensity function
$$
S\left(\vartheta ,t,\omega \right)=S_*+\vartheta _T  \int_{0
}^{t}h\left(t-s\right)\;{\rm d}X_s.
$$
{\sl Contiguous} means, that the likelihood ratio is {\sl asymptotically non degenerate}.  The function
$h\left(\cdot \right)$ is supposed to be known, bounded and
$$
h\left(\cdot \right)\in {\cal L}^1_+\left(\RR_+\right)=\left\{f\left(\cdot
\right)\geq 0\; : \;\int_{0}^{\infty } f\left(t\right)\;{\rm d}t<\infty \right\} .
$$

 To have contiguous alternatives we choose, as usual in regular problems,
$\vartheta_T =u/\sqrt{T}$, i.e.,
$$
S\left(u ,t,\omega \right)=S_*+ \frac{u}{\sqrt{T}} \int_{0
}^{t}h\left(t-s\right)\;{\rm d}X_s,\quad \quad u\geq 0.
$$
 Note that for any $h\left(\cdot \right)\in {\cal
L}^1_+\left(\RR_+\right)$ and any $u\leq K$ for sufficiently large $T$ the
condition \eqref{4} is fulfilled for the corresponding function $g\left(\cdot
\right)= uT^{-1/2}h\left(\cdot \right) $. This leads us to the following one
sided hypotheses testing problem:
\begin{eqnarray*}
&&{\scr H}_0:\qquad u=0,\qquad \left(\text{Poisson   process}\right)\\
&&{\scr H}_1:\qquad u>0,\qquad \left(\text{self-exciting process}\right).
\end{eqnarray*}
 This model corresponds to {\sl small self-exciting perturbations}
of the Poisson process of intensity $S_*$.

Note that as we use the LAN approach, we study the behavior of the tests statistics under hypothesis only (Poisson
process with constant intensity) and do not use the stationarity of the self-exciting processes under
alternatives. The limit of the  power function is obtained using LAN and Le
Cam's Third Lemma.

  Let us denote
$$
\Delta _\zT\left(X^T\right)=\frac{1}{S_*\sqrt{T\,{\rm I}_h^*}}\int_{0}^{T}\int_{0
}^{t-}h\left(t-s\right) \;{\rm d}X_s\;\left[{\rm d}X_t-S_*\;{\rm d}t\right].
$$
Here
$$
\int_{0}^{t-}h\left(t-s\right) \;{\rm d}X_s=\sum_{t_i<t}^{}h\left(t-t_i\right)
$$
(limit from the left of the integral, i.e., the term with $s_i=t$ is excluded)
 and
$$
{\rm I}_h^*=\int_{0}^{\infty }h\left(t\right)^2\;{\rm d}t+S_*\;\left(\int_{0}^{\infty
}h\left(t\right)\;{\rm d}t\right)^2
$$
is the Fisher information of the problem.  Throughout this paper we denote by
$z_\varepsilon $ the $1-\varepsilon $ quantile of the Gaussian law
$\ca{N}\left(0,1\right)$.
\begin{theorem}
\label{T1}
Let $h\left(\cdot \right)\in {\cal L}_+^1\left(\RR_+\right) $ and
bounded. Then the test
$$
\hat\phi_\zT\left(X^T\right)=\chi\zs{\left\{\Delta
_\zT\left(X^T\right)>z_\varepsilon \right\}}
$$
is locally asymptotically uniformly most powerful in the class ${\scr
 K}_\varepsilon $ and for any $u>0$ its
power function
\begin{equation}
\label{5}
\beta_\zT\left(u,\hat\phi_\zT\right)\longrightarrow
\hat\beta\left(u\right)=
\Pb\left\{ \zeta >z_\varepsilon -u\;\sqrt{{\rm I}_h^* }\right\},
\end{equation}
where $\zeta \sim \ca{N}\left(0,1\right)$.
\end{theorem}

\bigskip
\noindent
{\it Proof.} First note that the family of measures $\left\{\Pb_\vartheta
^{\left(T\right)}, \vartheta >0\right\}$ under hypothesis ${\scr H}_0 $ is LAN
at the point $\vartheta =0$, i.e., the random function
$Z_T\left(u\right)=L\left(u/\sqrt{T}, X^T\right)$ admits the representation
(see Kutoyants (1984), Theorem 4.5.3)
\begin{eqnarray*}
Z_T\left(u\right)&=&\exp\left\{
\int_{0}^{T}\ln \left(1+\frac{u}{S_*\sqrt{T}}\int_{0}^{t-}h\left(t-s\right){\rm
d}X_s \right)\;{\rm
d}X_t-
 \right.\\
&&\left.-
\frac{u}{\sqrt{T}}\int_{0}^{T} \int_{0}^{t}h\left(t-s\right)\;{\rm d}X_s\;{\rm
d}t\right\} =\\
&=&\exp\left\{u\sqrt{{\rm I}_h^*}\Delta_\zT\left(X^T\right)-
\frac{u^2}{2}\,{\rm I}_h^*+r_\zT\left(u,X^T\right)\right\}
\end{eqnarray*}
where
\begin{align}
\label{6}
\ca{L}_0\left\{\Delta_\zT\left(X^T\right) \right\}\Longrightarrow
\ca{N}\left(0,1\right)
\end{align}
and $r_\zT\left(u_\zT,X^T\right)\rightarrow 0 $ for any bounded sequence
$\left\{u_\zT\right\}$.

To verify \eqref{6} we check the following two conditions:
\begin{itemize}
\item {\sl Lindeberg  condition} for stochastic integral:
for any $\delta >0$
\begin{align*}
\lim_{T\rightarrow \infty }\frac{1}{T}\;\Ex_0\int_{0}^{T}
H_t^2\;\chi\zs{\left\{\left| H_t\right|>\delta \sqrt{T}\right\}}\;{\rm d}t=0,
\end{align*}
\item the {\sl law of large numbers}:
\begin{align}
\label{7}
\Pb_0-\lim_{T\rightarrow \infty }\frac{1}{S_*T}\int_{0}^{T}
H_t^2\;{\rm d}t={\rm I}_h^*.
\end{align}
\end{itemize}

Here we denoted
$$
H_t=\int_{0}^{t-} h\left(t-s\right)\;{\rm d}X_s.
$$
By these conditions the stochastic integral $\Delta _T\left(X^T\right)$ is
asymptotically normal. The proof of the corresponding central limit theorem
can be found, say,  in
Kutoyants (1984), Theorem 4.5.4  (of course, this theorem is a particular case of
general CLT for martingales).

To check these conditions we introduce an independent Poisson process
$\left\{X_t,\;t\leq 0\right\}$ of intensity $S_*$ and  replace $H_t $ by
$$
H_t^*=\int_{-\infty }^{t-} h\left(t-s\right)\;{\rm d}X_s.
$$
It is easy to see that for the process $H_t^*,t\geq 0$ we have
$$
\Ex_0\,H_t^*=S_*\;\int_{0 }^{\infty} h\left(v\right)\;{\rm d}v
$$
and
$$
\Ex_0\,\left(\left[H_t^*-\Ex_0\,H_t^*\right]\left[H_s^*-\Ex_0\,H_s^*\right]\right)
=S_* \;\int_{\max\left(0,s-t\right)}^{\infty } h\left(v+t-s\right)\;h\left(v\right)\;{\rm d}v.
$$
Note as well that
\begin{align*}
&\Pb_0\left\{ \frac{1}{\sqrt{T}}\int_{0}^{T}\left[H_t^*-H_t\right]\;
\left[{\rm d}X_t-S_*\,{\rm d}t\right]>\nu  \right\}\leq \\
&\quad \qquad \leq \frac{S_*}{T\nu ^2}\int_{0}^{T}
\Ex_0 \left(\int_{-\infty }^{0} h\left(t-s\right)\;{\rm d}X_s\right)^2{\rm
d}t=\\
&\quad \qquad =\frac{S_*^2}{T\nu ^2}\int_{0}^{T} \left[
\int_{t }^{\infty } h\left(v\right)^2\;{\rm d}v   +
S_*\left(\int_{t }^{\infty } h\left(v\right)\;{\rm d}v\right)^2
\right]\;{\rm d}t\longrightarrow 0
\end{align*}
as $T\rightarrow \infty $.

 Now the process $H_t^*,t\geq 0$ is second order  stationary and
$$
\Ex_0 \left(H_t^*\right)^2=S_*\int_{0}^{\infty }h\left(t\right)^2\;{\rm d}t
+S_*^2\left(\int_{0}^{\infty }h\left(t\right)\;{\rm d}t\right)^2= \Ex_0
\left(H_0^*\right)^2 <\infty .
$$
Hence
\begin{align*}
\Ex_0 \left(H_t^{*2}\;\chi\zs{\left\{\left| H_t^*\right|>\delta
\sqrt{T}\right\}}\right)=\Ex_0 \left(H_0^{*2}\;\chi\zs{\left\{\left|
H_0^*\right|>\delta \sqrt{T}\right\}} \right)\longrightarrow 0
\end{align*}
as $T\rightarrow \infty $ and
$$
\lim_{T\rightarrow \infty }\frac{1}{T}\int_{0}^{T}\Ex_0
\left(H_t^{*2}\;\chi\zs{\left\{\left| H_t^*\right|>\delta
\sqrt{T}\right\}}\right)\;{\rm d}t=0.
$$

The law of large numbers \eqref{7} will follow from the convergence:

\begin{align*}
M_T&=\Ex_0\left(\frac{1}{T}\int_{0}^{T}H_t^{*2}\;{\rm
d}t-\Ex_0\left(H_0^*\right)^2\right)^2 =\\
&=\frac{1}{T^2}\int_{0}^{T}\int_{0}^{T}
\Ex_0\left(H_t^{*2}-\Ex_0\left(H_0^*\right)^2
\right)\left(H_s^{*2}-\Ex_0\left(H_0^*\right)^2 \right)\;{\rm d}t\;{\rm
d}s\longrightarrow 0.
\end{align*}

To prove it we need the following elementary result.
\begin{lemma}
\label{L1}
Let $X=\left\{X_t,\;t\in A\right\}$ be a Poisson process of constant intensity
$S>0$ on $A\subset\RR$, and let $\ds f\left(\cdot \right),g\left(\cdot
\right)\in {\cal L}^k\left(A\right)=\left\{f\left(\cdot \right)\; : \;\int_{A}
\left|f(t)\right|^k\;{\rm d}t<\infty \right\}$, $k=1,\ldots,4$. Then
\begin{align*}
&\mathop{\bf Cov}
\left(\left(\int_{A}f\left(v\right)\;{\rm d}X_v\right)^2,
\left(\int_{A}g\left(v\right)\;{\rm d}X_v\right)^2\right)=\\
&\ =
4\int_{A}f(v)\,S\;{\rm d}v\,\int_{A}g(v)\,S\;{\rm
d}v\,\int_{A}f(v)\,g(v)\,S\;{\rm d}v+\\
&\quad+2\left(\int_{A}f(v)\,g(v)\,S\;{\rm
d}v\right)^2+\int_{A}f^2(v)\,g^2(v)\,S\;{\rm d}v+\\
&\quad+2\int_{A}f(v)\,S\;{\rm d}v\,\int_{A}f(v)\,g^2(v)\,S\;{\rm
d}v+2\int_{A}g(v)\,S\;{\rm d}v\,\int_{A}f^2(v)\,g(v)\,S\;{\rm d}v.
\end{align*}
\end{lemma}

\bigskip
\noindent
{\it Proof.} Using well-known properties of the Poisson processes (see, e.g.,
Kutoyants (1998), Lemma 1.1), we obtain the moment generating function
\begin{align*}
\phi\left(\lambda ,\mu \right)&=\Ex_0\exp\left\{\lambda
\int_{A}f\left(v\right)\;{\rm d}X_v+ \mu \int_{A}g\left(v\right)\;{\rm
d}X_v\right\} =\\
&=\exp\left\{\int_{A}\left(e^{\lambda f\left(v\right)+\mu
g\left(v\right)}-1 \right)\;S\;{\rm d}v\right\}.
\end{align*}
Remind that
\begin{align*}
\mathop{\bf Cov}
&\left(\left(\int_{A}f\left(v\right)\;{\rm d}X_v\right)^2,
\left(\int_{A}g\left(v\right)\;{\rm d}X_v\right)^2\right)=\\
&=\left.\frac{\partial ^4\phi\left(\lambda ,\mu \right) }{\partial\lambda
^2\:\partial\mu  ^2 }\right|_{\lambda =0,\mu =0}-
\left.\frac{\partial ^2\phi\left(\lambda ,0 \right) }{\partial\lambda
^2}\right|_{\lambda =0}
\left.\frac{\partial ^2\phi\left(0,\mu \right) }{\partial\mu
^2}\right|_{\mu=0}.
\end{align*}
Therefore the proof of the lemma follows from direct calculations.

\bigskip

Now we can write
\begin{align*}
R\left(t,s \right)&=\Ex_0\left(\left(H_{t
}^*\right)^2-\Ex_0\left(H_0^*\right)^2\right)
\left(\left(H_s^*\right)^2-\Ex_0\left(H_0^*\right)^2\right)=\\
&=4a^2K\left(t,s \right)+ 2K\left(t,s \right)^2+S_*\int_{-\infty }^{t\wedge s
}h\left(t-v\right)^2h\left(s-v\right)^2{\rm d}v  +\\
&\quad +
2aS_*\int_{-\infty }^{t\wedge s
}\left[h\left(t-v\right)^2h\left(s-v\right)+h\left(t-v\right)h\left(s-v\right)^2
\right]{\rm d}v
\end{align*}
where we put
$$
a=S_*\int_{0}^{\infty }h\left(y\right)\;{\rm d}y
$$
and (for $\tau =t-s$)
$$
K\left(t,s\right)=S_*\int_{\left|\tau\right|  }^{\infty
}h\left(y\right)\;h\left(y-\left|\tau\right| \right)\;{\rm d}y=K\left(\tau \right).
$$
Further, as the function $h\left(\cdot \right)$ is  bounded, we
have the estimate
$$
R\left(t,s  \right)\leq C\;K\left(\tau \right).
$$
 Hence
\begin{align*}
M_T=\frac{1}{T^2}\int_{0}^{T}\int_{0}^{T}R\left(t,s \right)\;{\rm d}t\;{\rm d}s
&\leq \frac{C}{T^2}\int_{0}^{T}\int_{0}^{T}K\left(t,s\right)\;{\rm d}t\;{\rm
d}s\leq \\
&\leq \frac{C}{T}\int_{-T}^{T}K\left(\tau \right)\;{\rm d}\tau .
\end{align*}
For the function $K\left(\cdot \right)$ we have
\begin{align*}
\int_{-T}^{T}K\left(\tau \right)\;{\rm d}\tau=S_*\int_{-T}^{T}
\int_{\left|\tau\right|}^{\infty}h\left(y\right)
\;h\left(y-\left|\tau\right|\right)\;{\rm d}y\;{\rm d}\tau
\leq 2\,S_*\left(\int_{0}^{\infty}h\left(y\right)\;{\rm d}y\right)^2.
\end{align*}

Hence $M_T\rightarrow 0$ and we have the law of large numbers \eqref{7}.

\bigskip

The property $\hat \phi_T\left(\cdot \right)\in {\scr K}_\varepsilon $ follows from the mentioned above asymptotic
normality of the statistic $\Delta_\zT\left(X^T\right) $.

Note as well that the convergence \eqref{5} follows from
$$
\ca{L}_u\left\{\Delta_\zT\left(X^T\right) \right\}\Longrightarrow
\ca{N}\left(u\;\sqrt{{\rm I}_h^*},\;1\right)
$$
(see the Third Lemma of Le Cam (van der Vaart (1998), p. 90)).

The asymptotic optimality of the test follows as well from the general theory
(see, e.g., Le Cam (1956) or Roussas (1972)), because if we replace $
{\scr H}_1$ by any
simple alternative $ {\scr H}_*:\; u=u_*$, then the test
$$
\bar \phi_T\left(X^T\right)=\chi\zs{\left\{L\left(u_*/\sqrt{T},
X^T\right)>b_\varepsilon \right\}}
$$
is the most powerful. Here
$$
b_\varepsilon =\exp\left\{u_*z_\varepsilon\sqrt{{\rm I}_h^*}-\frac{1}{2
}u_*^2\;{\rm I}_h^*\right\}\;\left(1+o\left(1\right)\right).
$$
It is easy to see that $\bar \phi_T\left(\cdot \right)\in  {\scr
K}_\varepsilon$ and the power function
$$
\beta_T \left(u_*,\bar \phi_T\right)\rightarrow \hat \beta \left(u_*\right).
$$
Therefore the test $\hat \phi_T\left(\cdot \right)$ is asymptotically as good
as the likelihood ratio test for any simple alternative.

\bigskip
\noindent
{\it Remark 1.} Note that the statistic $\Delta _T\left(X^T\right)$ can be
written as  follows
$$
\Delta _T\left(X^T\right)=\frac{1}{S^*\sqrt{T\,{\rm I}_h^*}}\sum_{0\leq t_j\leq T
}\sum_{t_i< t_j}
h\left(t_j-t_i\right)-\frac{1}{\sqrt{T\,{\rm I}_h^*}}\sum_{0\leq t_j\leq T
}\int_{0}^{T-t_j}h\left(v\right) \;{\rm d}v,
$$
where $t_i$ are the events of the observed process.

\bigskip
\noindent
{\it Remark 2.} By a similar way we can consider the problem of contiguous
hypotheses testing when under the hypothesis ${\scr H}_0$ the observed process
is self-exciting too. For example, let $h\left(\vartheta ,t\right)\geq 0,
t\geq 0$ be a smooth function of $\vartheta \in \Theta $, such that for all
$\vartheta \in \Theta $ the condition
$$
\int_{0}^{\infty }h\left(\vartheta ,t\right)\;{\rm d}t<1
$$
holds. Then  with the help of this function we introduce a family of
self-exciting processes with intensity functions
$$
S\left(\vartheta ,t,\omega \right)=S_*+\int_{-\infty }^{t}h\left(\vartheta
,t-s\right) \;{\rm d}X_s.
$$
Remind that these are stationary processes.

Now we can test the hypotheses
\begin{eqnarray*}
&&{\scr H}_0:\qquad \qquad \vartheta =\vartheta _0,\\
&&{\scr H}_1:\qquad \qquad \vartheta >\vartheta _0
\end{eqnarray*}
by the observations $X^T=\left\{X_t,0\leq t\leq T\right\}$.  Suppose as well
 that the function $h\left(\vartheta ,\cdot \right)$ is two
 times continuously differentiable on $\vartheta $ at the point $\vartheta
 =\vartheta _0$ and the derivatives $\dot h\left(\vartheta ,\cdot
 \right),\;\ddot h\left(\vartheta ,\cdot \right)$ satisfy the suitable
 conditions of integrability.   Let us denote
$$
\xi_t\left(\vartheta \right)=\int_{0 }^{t-}h\left(\vartheta
,t-s\right)\;{\rm d}X_s,\quad  \dot \xi_t\left(\vartheta \right)=\int_{0
}^{t-}
\frac{\partial h\left(\vartheta
,t-s\right)}{\partial \vartheta }\;{\rm d}X_s,
$$
and put
$$
\Delta _T\left(\vartheta _0,X^T\right)=\frac{1}{\sqrt{T}}\int_{0}^{T}\frac{\dot
\xi_t\left(\vartheta_0\right)}{ S_*+\xi_t\left(\vartheta_0\right)}
\;\left[{\rm d}X_t- S_*{\rm d}t-\xi_t\left(\vartheta _0\right)\;{\rm d}t\right].
$$

Then it can be easily shown that the test
$$
\hat \phi_T\left(X^T\right)=\chi\zs{\left\{\Delta _T\left(\vartheta
_0,X^T\right)> c_\varepsilon  \right\}}
$$
where $c_\varepsilon=z_\varepsilon \sqrt{{\rm I}_h\left(\vartheta _0\right)} $
is chosen from the condition $\hat \phi_T\in {\scr
K}_\varepsilon $ is locally asymptotically uniformly most powerful in the
class ${\scr K}_\varepsilon $. Here ${\rm I}_h\left(\vartheta _0\right)$ is the
Fisher information
$$
{\rm I}_h\left(\vartheta _0\right)=\Ex_{\vartheta
_0}\left(\frac{\dot{\xi}\left(\vartheta _0\right)^2 }{S_*+\xi\left(\vartheta
_0\right)}\right),
$$
where $\dot{\xi}\left(\vartheta _0\right) $ and $\xi\left(\vartheta
_0\right) $ are {\sl stationary random variables} related to the limit
distribution of the vector
$ \dot\xi_t\left(\vartheta_0\right), \xi_t\left(\vartheta_0\right)$.

\section{Testing of dependence}

Suppose that we have two sequences of events $0<t_1<t_2<\ldots <t_N<T$ and
$0<s_1<s_2<\ldots <s_M<T$ with corresponding counting processes
$X^T=\left\{X_t,0\leq t\leq T\right\}$ and $Y^T=\left\{Y_t,0\leq t\leq
T\right\}$. The first process is Poisson with constant known  intensity function
$S_X\left(t,\omega \right)=S_X>0 $ and the intensity function of the second
process can be written as
\begin{align*}
S_Y\left(t,\omega \right)=S_Y+\int_{-\infty }^{t}r\left(t-s\right)
\;{\rm d}X_s,
\end{align*}
where $r\left(\cdot \right)\in
{\cal L}^1\left(\RR_+\right)$.  Therefore, if $
r\left(t\right)\equiv 0$,
then the observed processes are standard (independent) Poisson processes of
intensities $S_X $ and $S_Y $ respectively (Hypothesis ${\scr H}_0$). For the
other values of $r\left(\cdot \right)$ we have
dependent point processes.

We suppose that the dependence between these two processes, if exists, is
weak, i.e., the function $
r\left(\cdot \right)$ is
sufficiently small  and we can apply the local approach. As before we suppose
that $r\left(t \right)=\vartheta _T\,h\left(t\right)$, where $h\left(\cdot
\right)\in {\cal L}^1\left(\RR_+\right)$ and $\vartheta
_T=u/\sqrt{T}\rightarrow 0$.
\begin{eqnarray*}
&&{\scr H}_0:\qquad u=0,\qquad \left({\rm independent \quad Poisson\quad
processes}\right)\\
&&{\scr H}_1:\qquad u>0,\qquad \left({\rm depending\quad  processes}\right).
\end{eqnarray*}

Introduce the statistic
\begin{align*}
\Delta_\zT&\left(X^T,Y^T\right)=\frac{1}{S_Y\sqrt{T\,{\rm I}_h}}\int_{0}^{T}\int_{0
}^{t-}h\left(t-s\right) \;{\rm d}X_s\;\left[{\rm d}Y_t-S_Y\;{\rm d}t\right]=\\
&=\frac{1}{S_Y\sqrt{T\,{\rm I}_h}}\sum_{0\leq s_j\leq T}\sum_{t_j<s_i}h
\left(s_j-t_i\right) -\frac{1}{\sqrt{T\,{\rm I}_h}}\sum_{0\leq t_j\leq T}
\int_{0}^{T-t_j}h\left(v\right)\;{\rm d}v
\end{align*}
where
$$
{\rm I}_h=\frac{S_X}{S_Y}\;\left(\int_{0}^{\infty }h\left(t\right)^2\;{\rm
d}t+S_X \left(\int_{0}^{\infty }h\left(t\right)\;{\rm
d}t\right)^2\right).
$$

\begin{proposition}
\label{P2}
Let $h\left(\cdot \right)\in {\cal L}^1_+\left(\RR_+\right)$ and bounded.
Then the test
$$
\hat\phi_\zT\left(X^T,Y^T\right)=\chi\zs{\left\{\Delta
_\zT\left(X^T,Y^T\right)>z_\varepsilon \right\}}
$$
is locally asymptotically uniformly most powerful in the class ${\scr
K}_\varepsilon $ and for any $u>0$ its power function
\begin{equation}
\label{7b}
\beta_\zT\left(u,\hat\phi_\zT\right)\longrightarrow
\hat\beta\left(u\right)=
\Pb\left\{ \zeta >z_\varepsilon -u\;\sqrt{{\rm I}_h }\right\},
\end{equation}
where $\zeta \sim \ca{N}\left(0,1\right)$.
\end{proposition}

\bigskip
\noindent
{\it Proof.} The proof is quite close to the given above proof of the
Theorem~1, and hence is omitted.

\bigskip
\noindent
{\it Remark 3.} The similar problem can be considered for the couple of
mutually exciting point processes with intensity functions
\begin{align*}
S_X\left(t,\omega \right)&=S_X+\int_{-\infty }^{t}r_{XY}\left(t-s\right)
\;{\rm d}Y_s,\\
S_Y\left(t,\omega \right)&=S_Y+\int_{-\infty }^{t}r_{YX}\left(t-s\right)
\;{\rm d}X_s,
\end{align*}
where $r_{XY}\left(\cdot \right),r_{YX}\left(\cdot \right)\in
{\cal L}^1\left(\RR_+\right)$.  Therefore, if $
r_{XY}\left(t\right)\equiv 0$ and $ r_{YX}\left(t\right)\equiv 0$,
then the observed processes are standard (independent) Poisson processes of
intensities $S_X >0$ and $S_Y>0 $ respectively (Hypothesis ${\scr
H}_0$). Under alternative there exists a weak dependence of these processes
through their intensity functions.

\section{One-sided nonparametric alternative}

In all considered above problems the alternatives are one-sided parametric. It
is possible to describe similar asymptotically uniformly most powerful tests
even in some nonparametric situations. Using the minimax approach we can
consider the least favorable model in the deriving of the upper bound on the
powers of all tests, but, of course, for special classes of intensities. This
approach sometimes is called semiparametric and the rate of convergence of
alternatives is $\sqrt{T}$.

As before, we suppose that under hypothesis ${\scr H}_0 $ the observed point
process $X^T=\left\{X_t,0\leq t\leq T\right\}$ is standard Poisson with known
intensity function $S\left(t\right)=S_*>0$ and under alternative $ {\scr H}_1$
it is self-exciting point process with intensity function
$$
S\left(t,\omega \right)=S_*+\int_{-\infty }^{t}g\left(t-s\right)\; {\rm
d}X_s,\quad  \quad 0\leq t\leq T
$$
where $g\left(\cdot \right)$ is now {\sl unknown} function. We suppose as well
that
\begin{equation}
\label{11}
\int_{0}^{\infty }g\left(t \right)\;{\rm d}t<1,
\end{equation}
hence the process $X^T$ is stationary. To describe the class of local
nonparametric alternatives we rewrite this intensity function as
$$
S\left(t,\omega \right)=S_*+\frac{1}{\sqrt{T}}\int_{-\infty
}^{t}u\left(t-s\right)\; {\rm d}X_s,\quad  \quad 0\leq t\leq T
$$
where the function $u\left(\cdot \right)$ is from the set ${\scr U}_{r}
$ defined below. Let us denote by ${\cal C}_+^b$ the set of nonnegative
functions bounded by the same constant and introduce  the set
$$
{\scr U}_r=\left\{u\left(\cdot \right)\in {\cal C}_+^b:\quad
\int_{0}^{\infty }u\left(t\right)\;{\rm d}t= r, \ \mathop{\rm
supp}u\left(\cdot\right)\text{ is bounded}\right\} .
$$
Note, that for any $r >0$ and $T> r ^2$ the condition \eqref{11}
is fulfilled.

Therefore, we consider the following hypotheses testing problem
\begin{eqnarray*}
&&{\scr H}_0:\quad\qquad u\left(\cdot \right)\equiv 0, \\
&&{\scr H}_1:\quad\qquad u\left(\cdot \right)\in {\scr U}_{r} ,\quad
r >0.
\end{eqnarray*}

The power  function of a test $\phi_\zT$ depends on the function
$u\left(\cdot \right)$ and we write it as
$$
\beta_T \left(u,\phi_\zT\right)=\Ex_u\;\phi_\zT \left(X^T\right).
$$
where $u=u\left(\cdot \right)\in {\scr U}_r$ with some $r >0$. We want to
apply an approach similar to the minimax one in the estimation theory. More
precisely, we seek to maximize the minimal power of test on the class ${\scr
U}_r$.  However, for any test $\phi_\zT\in {\scr K}_\varepsilon $ we have
$$
\inf_{u\left(\cdot \right)\in{\scr U}_r}\beta_T
\left(u,\phi_\zT\right)\leq \varepsilon
$$
 since for any $T>0$ we can take a function  from ${\scr U}_r$ equal
$0$ on $\left[0,T\right]$.
Hence we introduce the set
$$
{\scr U}_{r,N }=\bigl\{u\left(\cdot \right)\in{\scr U}_r:\quad
\mathop{\rm supp}u\left(\cdot \right)\subset\left[0,N\right] \bigr\},
$$
denote
$$
B_T \left(r ,N,\phi_\zT\right)= \inf_{u\left(\cdot
\right)\in {\scr U}_{r,N }}\beta_T \left(u,\phi_\zT\right)
$$
and give the following

\bigskip
\noindent
{\bf Definition 2.} A test $\phi^*_\zT\left(\cdot \right)$ is called {\sl
locally asymptotically uniformly most powerful in the class ${\scr
K}_\varepsilon$} if for any other test $\phi_\zT\left(\cdot \right)\in {\scr
K}_\varepsilon$ and any $K>0$ we have
$$
\lim_{N\rightarrow \infty }\lim_{T\rightarrow \infty }\inf_{0\leq r \leq
K }\left[B_T \left(r,N
,\phi_\zT^*\right)-B_T \left(r,N ,\phi_\zT\right) \right]\geq 0 .
$$

Let us introduce the decision function
$$
\hat \phi_\zT\left(X^T\right)=\chi\zs{\left\{\delta
_T\left(X^T\right)>z_\varepsilon \right\}},\qquad \qquad \delta
_T\left(X^T\right)=\frac{X_T-S_*T}{\sqrt{S_*T}}.
$$

\begin{theorem}
\label{T3}
The test $\hat \phi_\zT $ is locally asymptotically uniformly most powerful in
the class ${\scr K}_\varepsilon $ and for any $u\left(\cdot \right)\in {\scr
U}_r$ its power function
\begin{equation}
\label{12}
\beta_\zT\left(u,\hat\phi_\zT\right)\longrightarrow
\hat\beta\left(u\right)=\Pb\left\{ \zeta >z_\varepsilon -r
\;\sqrt{S_*}\right\},
\end{equation}
where $\zeta \sim \ca{N}\left(0,1\right)$.
\end{theorem}

\bigskip
\noindent
{\it Proof.} Let us fix a simple alternative $u\left(\cdot \right)\in {\scr
U}_r$, then the likelihood ratio
$
L_\zT\left(\frac{u\left(\cdot
\right)}{\sqrt{T}},X^T\right)=Z_T\left(u\left(\cdot \right)\right)
$
admits (under hypothesis ${\scr H}_0$) the representation (see the proof of
the theorem 1)
\begin{eqnarray*}
Z_T\left(u\left(\cdot \right)\right)&=&\exp\left\{ \int_{0}^{T}\ln
\left(1+\frac{1}{S_*\sqrt{T}}\int_{0}^{t-}u\left(t-s\right){\rm d}X_s
\right)\;{\rm d}X_t- \right.\\ &&\left.- \frac{1}{\sqrt{T}}\int_{0}^{T}
\int_{0}^{t}u\left(t-s\right)\;{\rm d}X_s\;{\rm d}t\right\} =\\
&=&\exp\left\{\Delta_\zT\left(u,X^T\right)- \frac{1}{2}\;{\rm I}\left(u\right)
+r_\zT\left(u,X^T\right)\right\}
\end{eqnarray*}
where
\begin{align*}
\Delta_\zT\left(u,X^T\right)&=\frac{1}{S_*\sqrt{T}}
\int_{0}^{T}\int_{0}^{t-}u\left(t-s\right){\rm d}X_s \;
\left[{\rm d}X_t-S_*\;{\rm d}t \right],\\
{\rm I}\left(u\right)&=\int_{0}^{\infty }u\left(t\right)^2\;{\rm
d}t+S_*\left(\int_{0}^{\infty }u\left(t\right)\;{\rm d}t\right)^2=
\int_{0}^{\infty }u\left(t\right)^2\;{\rm d}t+S_*\;r ^2
\end{align*}
and
$$
\ca{L}_0\left\{\Delta_\zT\left(u,X^T\right) \right\}\Longrightarrow
\ca{N}\left(0,{\rm I}\left(u\right)\right),\qquad
r_\zT\left(u,X^T\right)\rightarrow 0  .
$$
Moreover, these last two convergences are uniform on $u\left(\cdot \right)\in
{\cal U}_{r ,N}$, $0\leq r \leq K $ for any $K>0$. Hence the
likelihood ratio test
$$
\bar\phi_\zT\left(X^T\right)=\chi\zs{\left\{Z_T\left(u\left(\cdot
\right)\right)>d_\varepsilon  \right\}}
$$
with $d_\varepsilon=\exp\left\{z_\varepsilon \sqrt{{\rm I}\left(u\right) }-
{\rm I}\left(u\right)/2\right\} $
  is the most powerful in the class ${\scr K}_\varepsilon $ for any two
simple hypotheses and its power function
$$
\beta \left(u,\bar\phi_\zT\right)\longrightarrow \Pb\left\{\zeta
>z_\varepsilon - {\rm I}\left(u\right)^{1/2}\right\},\qquad \zeta \sim
{\cal N}\left(0,1\right).
$$
 It is easy to see that
$$
\inf_{u\left(\cdot \right)\in {\scr U}_{r,N} }{\rm I}\left(u\right)
=S_*\;r ^2+\frac{ r ^2}{N}
$$
because
$$
r^2 =\left(\int_{0}^{N}u\left(t\right)\;{\rm d}t\right)^2\leq N\
\int_{0}^{N}u\left(t\right)^2\;{\rm d}t
$$
with equality on the {\sl least favorable alternative}
$u^*\left(t\right)=\left(r /N\right)\chi\zs{\left\{0\leq t\leq N\right\}}$.

Hence
$$
\inf_{u\left(\cdot \right)\in {\scr U}_{r,N} }\Pb\left\{\zeta
>z_\varepsilon - {\rm I}\left(u\right)^{1/2}\right\}=\Pb\left\{\zeta
>z_\varepsilon - r\sqrt{S_*+N^{-1}} \right\}.
$$

Now we study the power function of the test $\hat\phi_\zT$. Let us denote
$$
U_t=\int_{0}^{t-}u\left(t-s\right)\,{\rm d}X_s,\qquad \pi _t=X_t-S_*t,
$$
then
$$
\Delta _T\left(u,X^T\right)=\frac{1}{S_*\sqrt{T}}\int_{0}^{T}U_t\,{\rm d}\pi _t,\qquad
\delta _\zT\left(X^T\right) =\frac{1}{\sqrt{S_*T}}\int_{0}^{T}\,{\rm d}\pi _t
$$
and
\begin{align*}
&\Ex_0\Delta _T\left(u,X^T\right)=0,\qquad
\Ex_0\Delta _T\left(u,X^T\right)^2={\rm I}\left(u\right),\quad
\Ex_0\delta_\zT\left(X^T\right)=0,\\
&\Ex_0\delta_\zT\left(X^T\right)^2=1\qquad \Ex_0 \left(\delta
_\zT\left(X^T\right)\Delta _T\left(u,X^T\right)\right) =r\sqrt{S_*} .
\end{align*}

Hence, under hypothesis ${\scr H}_0$, we have
$$
{\cal L}_0\left\{\Delta _T\left(u,X^T\right),\delta
_\zT\left(X^T\right)\right\}\Longrightarrow {\cal N}\left({\bf 0},{\bf R}\right)
$$
where $ {\bf R}$ is covariance matrix of the vector $\left(\Delta _T,\delta
_T\right) $ described above. Therefore $\hat\phi_\zT\in {\scr K}_\varepsilon$,
and using Le Cam's Third Lemma (van der Vaart (1998)) we obtain that under
alternative 
$u\left(\cdot \right)\in {\scr U}_r$
$$
\delta _\zT\left(X^T\right)\Longrightarrow {\cal N}\left(r \sqrt{S_*},1\right).
$$
For the power function we have
$$
\beta \left(u,\hat\phi_\zT\right)\longrightarrow \Pb\left\{\zeta
>z_\varepsilon - r\sqrt{S_*} \right\}.
$$
It can be shown that this convergence is uniform over $u\left(\cdot \right)\in
{\scr U}_{r,N}$, $0\leq r \leq K $ for any $K>0$ and this proves the
theorem.

\section{Simulations}

The main results (Theorems 1 and 2) of this work are {\sl asymptotic in
nature} and it is interesting to see the properties of the tests for the
moderate values of $T$. This can be done, say, by Monte-Carlo simulations.

\subsection{Parametric alternative}

To illustrate Theorem 1 we take $S_*=1$ and $h\left(t\right)=\frac{1}{2}{\rm
e}^{-t/2}$ (see Example 1). This yields
$$
S\left(u,t,\omega \right)=1+\frac{u}{2\sqrt{T}}\; \sum_{t_i\leq t}^{}{\rm
e}^{-\left(t-t_i\right)/2} , \qquad u\geq 0,\quad 0\leq t\leq T.
$$
In this case
$$
\Delta_\zT\left(X^T\right)=\frac{1}{\sqrt{5\,T}}\sum_{0\leq t_j\leq T}
\sum_{t_i<t_j}
{\rm e}^{-\left(t_j-t_i\right)/2}
-\frac{2}{\sqrt{5\,T}}\left(X_T-
\sum_{0\leq t_j\leq T}{\rm e}^{-\left(T-t_j\right)/2}\right)
$$
where $t_i$ are the events of the observed process, and the test
$\hat\phi_\zT^\varepsilon$ given by
$$
\hat\phi_\zT^\varepsilon=\hat\phi_{\zT}\left(X^T\right)=\chi\zs{\left\{\Delta
_\zT\left(X^T\right)>z_\varepsilon \right\}}
$$
is locally asymptotically uniformly most powerful in the class ${\scr
K}_\varepsilon$.

In Figure~1  we represent the size of the test $\hat\phi_\zT^{0.05}$ as a
function of $T\in[0,1000]$. This size is given by
$$
\alpha\left(T\right)=\Pb_0\left\{\Delta_T\left(X^T\right)>z_{0.05}\right\},
\qquad 1\leq T\leq 1000
$$
and is obtained by simulating $M=10^7$ trajectories on $[0,T]$ of Poisson
process of constant intensity $S\left(t,\omega \right)=1$ and calculating
empirical frequency of accepting the alternative hypothesis.

In Figure~2 we represent the power function of the test
$\hat\phi_\zT^{0.05}$ given by
$$
\beta_\zT\left(u,\hat\phi_\zT^{0.05}\right)=\Pb_u\left\{\Delta_T
\left(X^T\right)>z_{0.05}\right\},\qquad 0\leq u\leq 5
$$
for $T=100$, $300$ and $1000$, as well as the limiting (Gaussian) power
function given by
$$
\hat\beta\left(u\right)=\Pb\left\{\zeta
>z_{0.05}-u\,\sqrt{5}/2\right\}=\frac{1}{\sqrt{2\pi }}\int_{z_{0.05}
-u\,\sqrt{5}/2 }^{\infty } e^{-v^2/2}\;{\rm d}v,\quad 0\leq u\leq 5.
$$
The function $\beta_\zT$ is obtained by simulating (for each value of $u$)
$M=10^6$ trajectories on $[0,T]$ of self-exciting process of intensity
$S\left(u,t,\omega\right)$ and calculating empirical frequency of accepting
the alternative hypothesis.

\twopic{figure1.eps}{Fig.~1: Test size}{figure2.eps}{Fig.~2: Test power}

Now let us consider the $\tilde\phi_\zT^\varepsilon$ given by
$$
\tilde\phi_\zT^\varepsilon=\tilde\phi_{\zT}\left(X^T\right)=\chi
\zs{\left\{\Delta_\zT\left(X^T\right)>z\right\}}
$$
where the threshold $z$ is chosen so that this test is of exact size
$\varepsilon$. The choice of this threshold $z$ as a function of
$\varepsilon\in[0,0.25]$ is shown in Figures~3 and~4  for $T=100$, $300$
and $1000$, as well as the Gaussian threshold $z_\varepsilon$. The values of
$z$ are obtained by simulating $M=10^7$ trajectories on $[0,T]$ of Poisson
process of constant intensity $S\left(t,\omega \right)=1$ and calculating
empirical $1-\varepsilon$ quantiles of $\Delta_\zT$.

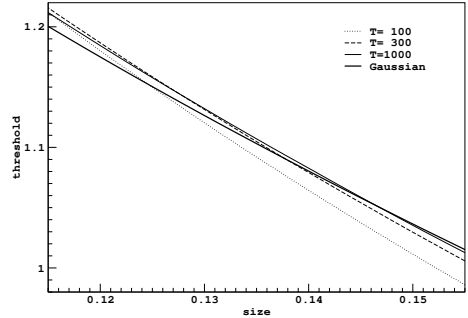
\twopic{figure3.eps}{Fig.~3: Threshold choice}{figure4.eps}{Fig.~4: Threshold
choice (zoom)}

For example to obtain test of exact size $0.05$ one needs take $z\simeq 1.78$
for $T=100$ ($z\simeq 1.74$ for $T=300$, $z\simeq 1.70$ for $T=1000$) against
$z_\varepsilon\simeq 1.64$ for Gaussian case.

\subsection{Nonparametric alternative}

To illustrate the nonparametric alternatives we take  intensity functions
corresponding to  $S_*=1$ and  $u\left(t\right)=\left(r/N
\right)\;\chi\zs{\left\{0\leq t\leq N\right\}}$, i.e.,
$$
S\left(t,\omega \right)=1+\frac{r }{N\sqrt{T}}
\sum_{t_i<t}^{}\chi\zs{\left\{ t-t_i\leq N\right\}} ,\qquad 0\leq t\leq T,
$$
where $t_i$ are the events of the observed process.  This choice of
$u\left(\cdot \right)$ allows us to compare the power function of our locally
asymptotically uniformly most powerful test
$$
\hat \phi_\zT^\varepsilon \left(X^T\right)=\chi\zs{\left\{ X_T> z_\varepsilon
\sqrt{T} +T\right\}}
$$
with the asymptotic power
$$
\beta \left(r \right)=\frac{1}{\sqrt{2\pi }}\int_{z_\varepsilon -r
}^{\infty } e^{-v^2/2}\;{\rm d}v
$$
of Neyman-Pearson test for the least favorable alternatives.

Note that under ${\scr H}_0$, $X_T$ is Poisson random variable with parameter
$T$, therefore the size of the test $\hat \phi_T^\varepsilon $, as well as the
threshold giving a test of exact size $\varepsilon$, can be calculated
directly (without resort to Monte-Carlo simulations).

We represent the power function of the test $\hat \phi_\zT^{0.05}$ given by
$$
\beta_\zT\left(u,\hat\phi_T\right)=\Pb_u\left\{X_T> z_{0.05}
\sqrt{T} +T \right\},\qquad 0\leq r\leq 5
$$
for $T=100$, $300$ and $1000$ as well as the limiting (Gaussian) function
$\beta \left(r \right)$, $0\leq r\leq 5$. In Figures~5 and~6 we take $N=5$ and
in $N=50$ respectively. The function $\beta_\zT$ is obtained by simulating
(for each value of $r$ and $N$) $M=10^6$ trajectories on $[0,T]$ of
self-exciting process of intensity $S\left(t,\omega\right)$ and calculating
empirical frequency of accepting the alternative hypothesis.

We see that if $1\ll N\ll T$, then the power function converge to the limiting
function (for example, if $N=50$ and $T=1000$, the power function almost
coincides with the limiting one).  If $N$ and $T$ are of the same order (for
example, if $N=50$ and $T=100$) then the power function of the test can be
essentially smaller. This example confirms the importance of use of functions
with bounded support and of the order of limits in Definition~2.

\twopic{figure5.eps}{Fig.~5: Test power ($N=5$)}{figure6.eps}{Fig.~6: Test
power ($N=50$)}

\section{Discussions}

The constructed tests are asymptotically optimal for parametric (Section 3)
and nonparametric (Section 5) alternatives. It seems that these  are just the
first results in this field and it is interesting to develope the
construction of the asymptotically optimal tests for  wider classes of
alternatives. Particularly, it is intersting to study  {\sl smooth
alternatives} like 
$$
{\scr H}_1:\qquad \int_{0}^{\infty
}u^{\left(k\right)}\left(t\right)^2{\rm d}t >r, 
$$
where $r>0$. Note that the test $\hat\phi_T$ is no more uniformly consistent
in this situation.

\bigskip

\noindent{\bf Acknowledgments}

\noindent We are grateful to anonymous referee and associate editor for useful comments.

\noindent
Yury A. Kutoyants, Laboratoire de Statistique et Processus,
Universit\'e du Maine, 72085 Le Mans Cedex 9, France, kutoyants@univ-lemans.fr

\end{document}